\documentclass{ifacconf}

\usepackage{graphicx}      
\usepackage{natbib}        
\usepackage{color}
\usepackage{textcomp}
\usepackage[mathcal]{euscript}
\usepackage{mathrsfs}
\usepackage{amsfonts,amssymb,amsbsy}
\usepackage{graphicx,empheq}
\usepackage{nomencl}
\usepackage{lipsum}
\usepackage{lmodern}
\newcommand{\mc}{\mathcal}

\renewcommand{\theenumii}{\roman{enumii}}
\renewcommand{\theenumiii}{\alph{enumiii}}
\renewcommand{\thepage}{}

\def\e1{{\varepsilon_{11}}}
\def\b1{{\beta_{11}}}
\def\bp3{{\beta_{33}}}
\def\ep3{{\varepsilon_{33}}}

\def\Ltwo{{\mathbb L}^2 }

\MHInternalSyntaxOn
\providecommand*\phantomword[3][c]{%
\mathchoice
{\MT_phantom_word:NNnn #1\displaystyle {#2}{#3}}%
{\MT_phantom_word:NNnn #1\textstyle {#2}{#3}}%
{\MT_phantom_word:NNnn #1\scriptstyle {#2}{#3}}%
{\MT_phantom_word:NNnn #1\scriptscriptstyle {#2}{#3}}%
}
\def\MT_phantom_word:NNnn #1#2#3#4{%
\@begin@tempboxa\hbox{$\m@th#2#4$}%
\setlength\@tempdima{\widthof{$\m@th#2#3$}}%
\hbox{\hb@xt@\@tempdima{\csname bm@#1\endcsname}}%
\@end@tempboxa}
\MHInternalSyntaxOff

\begin{document}
\begin{frontmatter}

\title{Exponential stabilization of a smart piezoelectric composite beam with only one
boundary controller \thanksref{footnoteinfo}}

\thanks[footnoteinfo]{ This  research is supported by the Western Kentucky University startup grant. }

\author[First]{Ahmet \"{O}zkan \"{O}zer}

\address[First]{Department of Mathematics, Western Kentucky University,\\
   Bowling Green, KY 42101, USA (E-mail:ozkan.ozer@wku.edu).}

\begin{abstract}                
Layered smart composite beams involving a piezoelectric layer are traditionally actuated by a voltage source by the extension mechanism. In this paper, we consider only the bending and shear of a cantilevered piezoelectric smart composite beam  modeled by the Mead-Marcus sandwich beam assumptions.  Uniform exponential stabilitization with only one boundary state feedback controller, simultaneously controlling both bending moment and shear, is proved by using a spectral multiplier approach. The state feedback controller slightly differs from the classical counterparts by a non-trivial compact and nonnegative integral operator. This is due to the strong coupling of the charge equation with the stretching and bending equations. For simulations, the so-called filtered semi-discrete finite difference scheme  is adopted.
\end{abstract}

\begin{keyword}
Piezoelectric smart composite, smart sandwich beam, boundary feedback stabilization, electrostatic, Mead-Marcus sandwich beam.
\end{keyword}

\end{frontmatter}

\section{Introduction}
\begin{figure}[h!tb]
\centering
\includegraphics[height=1.5in]{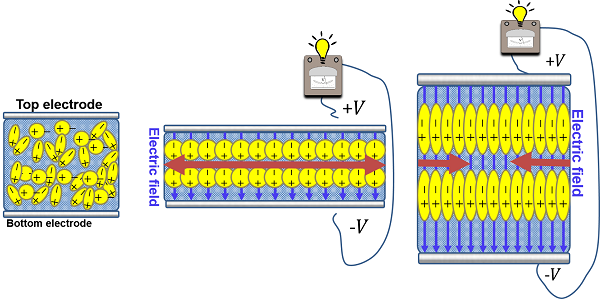}
\includegraphics[height=1.7in]{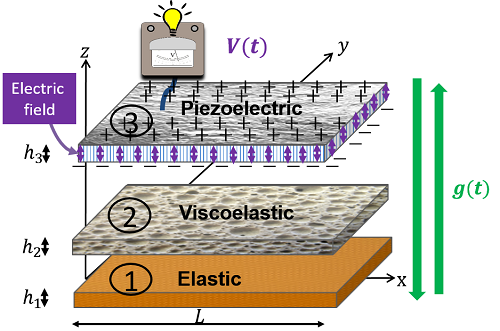}
\caption{\scriptsize (a) A piezoelectric beam extents or shrinks by supplying voltage to its electrodes since the charges separate and line up in the vertical direction.  (b) A voltage-actuated piezoelectric smart composite of length $L$  with thicknesses $h_1, h_2, h_3$ for its layers { \textcircled{1},  \textcircled{2}, \textcircled{3},} respectively.  Both voltage $V(t)$ and shear $g(t)$  controller control bending motions on the composite. In fact, it is the goal of the paper that the voltage controller $V(t)$ itself has the ability to control all bending and shear motions on the composite  in a few seconds. 
} \label{ACL}
 \end{figure}

A piezoelectric smart composite beam is a three-layer sandwich beam consisting of a stiff elastic layer, a complaint (viscoelastic) layer, and a piezoelectric layer, see Fig. \ref{ACL}.
The piezoelectric layer is also an elastic beam  with  electrodes at its top and bottom surfaces and connected to an external electric circuit.  As the electrodes are subjected to a voltage source, an electric field is created between the electrodes, and the piezoelectric beam  shrinks or extends. Therefore, the whole composite stretches and bends  (see Fig. \ref{ACL}). 

The modeling assumptions for smart piezoelectric models can be classified in two main categories: mechanical and electro-magnetic. The mechanical assumptions can be classified in two main categories, i.e. see \cite{Trindade}: either Mead-Marcus (M-M) type \cite{Baz} or  Rao-Nakra (R-N) type \cite{Baz,Lam}.  The M-M models only involve the transverse kinetic energy whereas the R-N models involve both longitudinal and transverse kinetic energies. Both types of models  reduce to the classical counterparts, see  \cite{OzerACC}, once the piezoelectric strain is taken to be zero. The electro-magnetic assumptions on the piezoelectric layer are either fully dynamic, quasi-static, or electrostatic, see \cite{O-M1,Ozkan1}. The electrostatic assumption completely discards electrical and magnetic-kinetic energies due to Maxwell's equations. It is still a standard assumption in the literature, see \cite{Smith}. The voltage control, actuating the piezoelectric layer, is simply  blended into models through the boundary conditions.

For the passive sandwich beam models (having no piezoelectric layer), the exact controllability of the M-M and R-N models are shown for the clamped and hinged models \cite{O-Hansen1,O-Hansen3}. The exponential stability in the existence of the passive damping term due to the shear of the middle layer is investigated for the  M-M model (\cite{A-H,W-G}) . The active boundary feedback stabilization of the classical R-N model is only investigated  for hinged (\cite{O-Hansen4}) and clamped-free (\cite{Wang}) boundary conditions.
The exponential stabilizability of the cantilevered fully dynamic or electrostatic M-M and R-N models has been open problems for more than a decade. Note that cantilevered boundary conditions are more physical than clamped or hinged boundary conditions. Recently,  the  exponential stability of the electrostatic  R-N model is shown  by using four feedback controllers \cite{OzerIEEE}, two for stretching motions of outer layers, and two for the bending motion. The exponential stability with only three controllers is recently shown by using a  spectral-theoretic approach  \cite{Y-W}, and by a higher order spectral  multipliers approach (\cite{OzerIEEE,Ozer17}). The fully dynamic R-N model is shown to be not stabilizable for many choices of material parameters by using   $B^*-$type feedback controllers  \cite{OzerIEEE}.  The charge-actuated electrostatic counterparts are also  shown to be exponentially stable in \cite{Ozer18a}.

To our knowledge, the exponential stabilizability for ``cantilevered" fully dynamic or electrostatic M-M model have never been studied in the literature. Denoting  stretching of the top and the bottom layers, bending of the composite, shear due to the middle layer, and the total induced charge accumulated at the piezoelectric layer by $v^1,v^3,w,\phi^2,p$ respectively, the equations of motion for the fully dynamic M-M model  is obtained in \cite{Ozer17} by a thorough variational approach as the following
\begin{eqnarray}
\nonumber  && \left\{ \begin{array}{ll}
  m \ddot w + A w_{xxxx} -  B_1\varsigma \gamma\beta h_2 h_3 \phi^2_x + \gamma\beta B_3p_{xxx}=0,&\\
  C\varsigma \phi^2 -\phi^2_{xx} + B_1w_{xxx} + B_2 p_{xx}=0,&\\
 \mu  h_3 \ddot p   - B_4 \beta p_{xx} + \gamma \beta h_2h_3 \varsigma B_2\phi^2&\\
 \quad\quad\quad\quad\quad -\gamma \beta B_3 w_{xxx}  = -V(t)\delta_L,&
 \end{array}\right.\\
  \label{perturbed-dumb} &&\left\{\begin{array}{ll}
 \left| w, w_x, \phi^2, p~\right|_{x=0},\quad \left| w_{xx}= \phi^2_x=p_x\right|_{x=L}=0, &\\
 \left| - A w_{xxx} + B_1\varsigma \gamma \beta h_2 h_3 \phi^2 - \gamma\beta B_3p_{xx}\right|_{x=L}=g(t).
\end{array}\right.
\end{eqnarray}
where $\delta_L=\delta(x-L)$ is the Dirac-Delta distribution at $x=L,$
$\phi^2 =\frac{1}{h_2}\left(-v^1+v^3\right) + \frac{H}{h_2}w_x,$ $H=\frac{h_1 + 2h_2+h_3}{2},$ and $h_i$ is the thickness of the ${\rm{i}}^{\rm th}-$ layer, and $\beta,\gamma, \mu>0$ are piezoelectric constants,  and \\$m, A, B_1, \ldots, B_4, C,\varsigma >0$ are functions for material  parameters of each layer. Moreover, $V(t)$ is the voltage controller actuating the piezo-layer, and $g(t)$ is actuating the transverse shear mechanism at the tip. The lack of stabilizability of this model for certain sub-classes of solutions is studied in \cite{Ozer17}.

Notice that if the electrostatic assumption is adopted, i.e. $ \mu h_3\ddot p\equiv 0,$ the model (\ref{perturbed-dumb}) reduces to
 \begin{eqnarray}
\nonumber  && \left\{\begin{array}{ll}
m \ddot w + {\tilde A} w_{xxxx} -  \beta \gamma h_2h_3 \varsigma{\tilde B}  \phi^2_x=-\frac{\gamma B_3}{B_4}V(t) (\delta_L)_x,& \\
\varsigma{\tilde C} \phi^2 -\phi^2_{xx} + {\tilde B} w_{xxx} =-\frac{B_2}{\beta B_4}V(t) \delta_L,&
\end{array}\right.\\
\label{abstractMM} &&\left\{ \begin{array}{ll}
  w(0)=w_x(0)=\phi^2(0)= w_{xx}(L)=\phi^2_x(L)=0,&\\
 - \tilde A w_{xxx}(L) + \beta \gamma h_2h_3 \varsigma{\tilde B} \phi^2(L)=g(t)
 \end{array}\right.
 \end{eqnarray}
where  the coefficients are $\tilde A=A-\frac{\gamma^2 \beta B_3^2}{B_4}>0, \tilde B=B_1-\frac{\gamma B_2 B_3}{B_4}>0, \tilde C=C+ \frac{\gamma h_2 h_3 B_2^2}{B_4},$ see \cite{OzerACC}.

 The case $g(t) \ne 0, V(t)=0$ corresponds to the standard (passive) M-M model, and its stabilizability is studied in \cite{Wang}. To our knowledge, the only stabilizability result for the electrostatic model ($V(t)\ne 0$, $g(t)=0$) is provided by \cite{Baz} where various PID-type feedback controllers are considered for the asymptotic stability of the system. These results  do not imply the exponential stability whatsoever. In fact, a shear-type of passive damping is also included in their models as the following:
 \begin{eqnarray}
\label{abstractMM-damped} \left\{ \begin{array}{ll}
m \ddot w + {\tilde A} w_{xxxx} -  \beta \gamma h_2h_3 \varsigma{\tilde B}  \phi^2_x=-\frac{\gamma B_3}{B_4}V(t) (\delta_L)_x, &\\
  \kappa\dot \phi^2 + \varsigma{\tilde C} \phi^2 -\phi^2_{xx} + {\tilde B} w_{xxx} =-\frac{B_2}{\beta B_4}V(t) \delta_L
  \end{array}\right.
 \end{eqnarray}
where $\kappa>0$ is the damping coefficient. It is proven in  \cite{W-G} that the damping term itself exponentially dissipates the energy of (\ref{abstractMM-damped}), even without the boundary feedback damping: $V(t)\equiv 0$. Hence, it is not clear whether $V(t)$ can be designed to exponentially dissipate the energy by itself.

In this paper, first we show that the the model (\ref{abstractMM}) is well-posed on an appropriate Hilbert space. Next, we prove that the overdetermined problem, with an extra measurement, has only the trivial solution by using spectral multipliers to ensure the strong stability.  Without  considering the shear-type of passive damping, i.e. $\kappa=0$ in (\ref{abstractMM-damped}), the exponential stability of the electrostatic M-M model is guaranteed by using only the $B^*-$type state feedback controller for $V(t)$. The proof  combines the a spectral multiplier  method and a frequency domain approach as in \cite{L-L}. Finally, the so-called filtered semi-discrete  Finite Differences is proposed first time to design the approximated stabilizing controller for a strongly coupled system. 

\section{Well-posedness}
Define the operator $(\varsigma \tilde C I- D_x^2)$ on the domain $H^2_L(0,L):=\{\psi\in H^2(0,L): \psi_x(0)=\psi_x(L)=0\}.$ Therefore, the operator $P_\varsigma=(\tilde C\varsigma I-D_x^2)^{-1}$ is defined by
\begin{eqnarray}
  \label{def} && P_\varsigma (f)(x) = \int_0^L g(x,z)f(z) dz,\quad {\rm with} \\
\nonumber &&g(x,z) = \left\{
            \begin{array}{ll}
           \frac{\cosh{[\sqrt{\tilde C\varsigma}(z-L)]} \sinh{[\sqrt{\tilde C\varsigma} x]}}{\sqrt{\tilde C\varsigma}\cosh{[\sqrt{\tilde C\varsigma}L]}} & x\le z \\
          \frac{\cosh{[\sqrt{\tilde C\varsigma}(x-L)]} \sinh{[\sqrt{\tilde C\varsigma} z]}}{\sqrt{\tilde C \varsigma}\cosh{[\sqrt{\tilde C\varsigma}L]}} & x\ge z.
    \end{array}
  \right.
\end{eqnarray}
It is well-known that $P_\varsigma$ is a compact  and non-negative  operator on $\Ltwo(0,L)$. We have the following result:

\begin{lem} \label{pxi} Let ${\rm Dom} (D_x^2)=H^2_L(0,L).$ Define the operator $J_\varsigma:=\varsigma \tilde C P_\varsigma-  I. $ Then, $J_\varsigma$ is continuous, self-adjoint and and non-positive on $\Ltwo(0,L).$ Moreover, for all $w\in {\rm Dom} (P_\varsigma), $
 $J_\varsigma w =  ~P_\varsigma D_x^2  =  (\varsigma \tilde C I- D_x^2)^{-1}D_x^2  w.$
\end{lem}
\noindent {\bf Proof:} Continuity and self-adjointness easily follow from the definition of $J_\varsigma.$ We first prove that $J_\varsigma$ is a non-positive operator. Let $u\in \Ltwo(0,L).$ Then $P_\varsigma u=(\varsigma \tilde C I- D_x^2)^{-1} u=s$ implies that $s\in {\rm Dom} (D_x^2)$ and $\varsigma \tilde C s- s_{xx}=u$
\begin{eqnarray}\nonumber \left<J_\varsigma u, u\right>_{\Ltwo(0,L)}&=& \left<(\varsigma \tilde C P_\varsigma-I) u, u\right>_{\Ltwo(0,L)}\\
\nonumber &=& -\varsigma \tilde C \|  s_{x} \|^2_{\Ltwo(0,L)} - \|  s_{xx}\|^2_{\Ltwo(0,L)}.
\end{eqnarray}
Let $J_\varsigma w=(\varsigma \tilde C P_\varsigma- I) w$ and $v:=P_\varsigma   w.$ Then $\varsigma \tilde C v- v_{xx}=w.$ By a simple rearrangement of the terms $$J_\varsigma w= (\varsigma \tilde C v- w)=(\varsigma \tilde C v-\varsigma \tilde C v+ v_{xx})=v_{xx}= P_\varsigma D_x^2 w. \square$$

By Lemma \ref{pxi}, (\ref{abstractMM}) can be simplified to
\begin{eqnarray}
\label{MM} \left\{\begin{array}{ll}
 m \ddot w + {\tilde A} w_{xxxx}+\gamma\beta \varsigma h_2 h_3 \tilde B^2 (J_\varsigma w_{x})_x=&\\
 \quad-\frac{\gamma}{B_4}\left[ \varsigma h_2 h_3  \tilde B B_2 (P_\varsigma \delta_L)_x + B_3(\delta_L)_x\right]V(t),\\
  w(0)=w_x(0)= w_{xx}(L)=0, &\\
  {\tilde A} w_{xxx}(L)+\gamma\beta \varsigma h_2 h_3 \tilde B^2 J_\varsigma w_{x}(L)=g(t).
 \end{array} \right.
    \end{eqnarray}
  This model fits in the form of the abstract Mead-Marcus beam model obtained in \cite{O-Hansen1}.

 Since our beam in nonclassical, we discard the mechanical controller; $g(t)\equiv 0.$ Define $\mc{H} =\mathrm V \times \mathrm H= H^2_L(0,L) \times \Ltwo(0,L).$
The  energy associated with (\ref{MM}) is
\begin{eqnarray}
\nonumber  \mathrm{E} =\frac{1}{2}\int_0^L \left\{ m |\dot w|^2  + \tilde A |w_{xx}|^2 -\gamma\beta \varsigma h_2 h_3 \tilde B^2 J_\varsigma w_{x} \bar w_x \right\} ~dx.
\end{eqnarray}
 This motivates the definition of the inner product on $\mc{H}:$
\begin{eqnarray}
\nonumber && \left<\left[ \begin{array}{l}
 u_1 \\
 u_2
 \end{array} \right], \left[ \begin{array}{l}
 v_1 \\
 v_2
 \end{array} \right]\right>_{\mc{H}}= \left<u_2, v_2\right>_{\mathrm H} + \left<u_1, v_1\right>_{\mathrm V}= \int_0^L \left\{m  u_2 { {\bar v}_2}\right.\\
 \nonumber  &&\left.+ \tilde A (u_2)_{xx} (\bar v_2)_{xx}  -\gamma\beta \varsigma h_2 h_3 \tilde B^2 (J_\varsigma (u_1)_{x}) ({\bar u}_1)_x  \right\}~dx.
 \end{eqnarray}
Define  the operator $\mc A: {\text{Dom}}(\mc A)\subset \mc H \to \mc H$
where $\mc A= \left[ {\begin{array}{*{20}c}
   0 & I \\
       \frac{-1}{m}\left(\tilde A D_x^4+\gamma\beta \varsigma h_2 h_3 \tilde B^2 D_x J_\varsigma D_x\right) &  0  \\
\end{array}} \right]$ with
\begin{eqnarray}
 \label{A-MM-newd}  \left.
\begin{array}{ll}
{\rm {Dom}}(\mc A) = \{ (z_1,z_2)\in \mc H, z_2\in H^2_L(0,L),~~\tilde A (z_1)_{xxx}&\\
+\gamma\beta \varsigma h_2 h_3 \tilde B^2  J_\varsigma (z_1)_x \in H^1(0,L), ~~ (z_1)_{xx}(L)=0 &\\
 \tilde A (z_1)_{xxx}(L)+\gamma\beta \varsigma h_2 h_3 \tilde B^2  J_\varsigma (z_1)(L)=0 \}. &
 \end{array} \right.
\end{eqnarray}
 Define also the control operator $\mc B \in \mathcal{L}(\mathbb{C} , {\rm Dom}(\mc A)')$ by
 \begin{eqnarray}
 \label{defb_000} \mc B=   \left[ \begin{array}{c} 0 \\ -\frac{\gamma}{mB_4}\left[ \varsigma h_2 h_3  \tilde B B_2(P_\varsigma \delta_L)_x + B_3(\delta_L)_x\right] \end{array} \right].
 \end{eqnarray}
 The dual operator $\mc B^*\in \mathcal{L}( \mc H, \mathbb C)$ is defined by
 $\mc B^*\Phi=\frac{\gamma}{m B_4}\left[ \varsigma h_2 h_3  \tilde B B_2(P_\varsigma (\Phi_2)_x(L))  + B_3(\Phi_2)_x(L)\right].$\\
Choosing the state $\Phi=[w,\dot w]^{\rm T},$ the control  system (\ref{MM})  with the voltage controller $V(t)$ can be put into the  state-space form
\begin{eqnarray}
\label{Semigroupp-mmm}
\dot \Phi = {\mc A}  \Phi +\mc B V(t) , \quad \Phi(x,0) =  \Phi ^0.
\end{eqnarray}
Since the piezoelectric smart beam model is similar to the classical counterpart with the electrostatic assumption, the following results are immediate from (\cite{Ozer17}):
\begin{thm}For fixed initial data and no applied forces, the solution $(w,p)\in \mc H$ of (\ref{perturbed-dumb}) converges to  the solution of $(w,p)\in \mc H$ in (\ref{MM}) as $\mu\to 0.$
\end{thm}
\begin{thm}\label{w-pff}
Let $T>0,$ and $V(t)\in \Ltwo(0,T).$ For any $\Phi^0 \in \mathrm{H},$ $\Phi\in C[[0,T]; \mc H]$ and there exists a positive constant $c_1(T)$
such that (\ref{Semigroupp-mmm}) satisfies
      \begin{eqnarray}\label{conc}\|\Phi (T) \|^2_{\mc{H}} &\le& c_1 (T)\left\{\|\Phi^0\|^2_{\mc H} + \|V\|^2_{\Ltwo(0,T)}\right\}.
      \end{eqnarray}
\end{thm}

\section{Uniform Stabilization}

For $k_1>0,$ we choose the following $B^*-$type feedback controller
\begin{eqnarray}\nonumber && V(t)=-\frac{k_1 m B_4}{\gamma}\mc B^*\Phi\\
 \label{feed}&&=- k_1\left[ \varsigma h_2 h_3  \tilde B B_2(P_\varsigma \dot w_x(L))  + B_3\dot w_x(L)\right].
\end{eqnarray}
The energy of the system is dissipative and it satisfies
\begin{eqnarray}
\nonumber&& \frac{dE(t)}{dt}=\gamma V(t) \left[ h_2 h_3 \varsigma \tilde B (P_\varsigma \dot w_x)(L) +  \frac{B_3}{B_4} \dot w_x(L)\right]\\
\nonumber &&=-k_1\left[ \left( h_2 h_3 \varsigma \tilde BB_2 P_\varsigma  +  B_3 I\right)\dot w_x(L) \right]^2\le 0
\end{eqnarray}
where $ h_2 h_3 \varsigma \tilde BB_2 P_\varsigma  +  B_3 I$ is a non-negative operator. 

Observe that $P_\varsigma \dot w_x(L)$  is a PID-type feedback, and it is the total piezoelectric effect  due to the coupling of the charge equation to shear and bending at the same time. By  Lemma \ref{pxi}, it can also be considered as
$P_\varsigma \dot w_x(L) = \frac{1}{\varsigma \tilde B \tilde C}(-\dot \phi^2+\tilde B \dot w_x)(L).$
This type of representation is helpful to design the controller numerically inSection \label{simul}.
 Therefore (\ref{feed}) reduces to
\begin{eqnarray*} V(t)=-k_1\left[\left( \frac{h_2h_3\tilde B B_2}{\tilde C}+B_3\right) \dot w_x(L) - \frac {h_2h_3 B_2}{\tilde C} \dot \phi^2(L)\right].
\end{eqnarray*}

Now consider the system (\ref{Semigroupp-mmm}) with the state feedback controller (\ref{feed}):
\begin{eqnarray}
\label{Semigroupp-f}\left\{
\begin{array}{ll}
\dot \varphi = \tilde A\Phi:= \left({\mc A} -\frac{k_1 m B_4}{\gamma}\mc B \mc B^*\right) \Phi, \quad \Phi(x,0) =  \Phi ^0.\quad
\end{array}\right.
\end{eqnarray}
\begin{thm} The operator $\tilde {\mc A}$ defined by (\ref{Semigroupp-f}) is dissipative in $\mc H.$ Moreover,
${\tilde {\mc A}}^{-−1}$ exists and is compact on $\mc H.$ Therefore, $\tilde {\mc A}$ generates a $C_0$-semigroup of contractions on $\mc H$ and the
spectrum $\sigma(\tilde {\mc A})$ consists of isolated eigenvalues only.
\end{thm}

{\bf Proof} Let $Y\in{\rm Dom}(\tilde {\mc A}).$ Then
\begin{eqnarray*}
\nonumber &&\left<\mc AY,Y\right>=\left. \left(-\tilde A (y_1)_{xxx} -\gamma\beta \varsigma h_2 h_3 \tilde B^2 D_x J_\varsigma  (y_1)\right) \bar y_2\right|_{x=0}^L\\
\nonumber &&\quad+ \left.\tilde A (y_1)_{xx} (\bar y_2)_x\right|_{x=0}^L\\
\nonumber &&\quad +\int_0^L\left[ -\tilde A (y_1)_{xx} (\bar y_2)_{xx}+ \gamma\beta \varsigma h_2 h_3 \tilde B^2 J_\varsigma (y_1)_x  (\bar y_2)_x \right.\\
\nonumber &&\quad\left.+\left(\tilde A (y_2)_{xx} (\bar y_1)_{xx}  -\gamma\beta \varsigma h_2 h_3 \tilde B^2 J_\varsigma (y_2)_{x}  (\bar y_1)_x \right)\right]~dx.
\end{eqnarray*}
Therefore,
\begin{small}
\begin{eqnarray}\label{eq7}{\rm Re}\left<\tilde {\mc A} Y,Y\right>=-\left[ \left( h_2 h_3 \varsigma \tilde BB_2 P_\varsigma  +  B_3 I\right)(y_2)_x(L) \right]^2\le 0.\quad~~
\end{eqnarray}
\end{small}
Therefore $\tilde {\mc A}$ is dissipative. If ${\tilde {\mc A}}^{-−1}$ exists, $\mc A$ must be densely defined in $\mc H.$ Therefore, $\tilde {\mc A}$ generates a $C_0$-semigroup of contractions on $\mc H.$ Next, we show that $0\in \sigma (\tilde {\mc A}),$ i.e. $0$ is not an eigenvalue. We solve the following problem:
\begin{eqnarray}
 \label{eig1-zero}\left\{\begin{array}{ll}
 {\tilde A} w_{xxxx}+\gamma\beta \varsigma h_2 h_3 \tilde B^2 (J_\varsigma w_{x})_x  = 0, &\\
 w(0)=w_x(0)= w_{xx}(L)= 0,&\\
 {\tilde A} w_{xxx}(L)+\gamma\beta \varsigma h_2 h_3 \tilde B^2 J_\varsigma w_{x}(L)=0.
 \end{array}\right.
\end{eqnarray}
 Let $J_\varsigma w_x:=u.$ By the definition of $J_\varsigma=(\varsigma \tilde C I -D_x^2)^{-1} D_x^2,$ (\ref{eig1-zero}) is re-written as
  \begin{eqnarray}
\nonumber \left\{ \begin{array}{ll}
 {\tilde A} w_{xxxx} -  \beta \gamma h_2h_3 \varsigma{\tilde B}  u_x=0, &\\
\varsigma{\tilde C} u -u_{xx} + {\tilde B} w_{xxx} =0,&\\
w(0)=w_x(0)=u(0)= w_{xx}(L)=u_x(L)=0,&\\
  \tilde A w_{xxx}(L) - \beta \gamma h_2h_3 \varsigma{\tilde B} u(L)=0.
\end{array}\right.
 \end{eqnarray}
By using the last boundary condition, we integrate the first equation and plug it in the $u-$equation to get
$\left(\xi+\frac{\beta \gamma h_2h_3 \varsigma{\tilde B}^2  }{\tilde A}\right)u-u_{xx}=0.$
Since $\xi +\frac{\beta \gamma h_2h_3 \varsigma{\tilde B}^2  }{\tilde A} >0,$ by the boundary conditions for $u,$ we obtain that $u\equiv 0.$ This implies that $w_{xxx}=0.$ By the boundary conditions $w\equiv 0.$
$ {\tilde A} w_{xxx} -  \beta \gamma h_2h_3 \varsigma{\tilde B}  u=0.$
 Thus, $0\in \sigma (\mc A),$ and $\mc A^{−-1}$ is compact on $\mc H.$  Hence the
spectrum $\sigma(\mc A)$ consists of isolated eigenvalues only. $\square$

\begin{thm}
\label{strong}
   The solutions $\Phi(t)$  for $t\in \mathbb{R}^+$ of the closed-loop system (\ref{Semigroupp-f})
     is strongly stable in $\mc H.$
\end{thm}

\noindent {\bf Proof:} If we can show that
there are no eigenvalues on the imaginary axis, or in other words, the set
\begin{eqnarray}\label{set} &&\left\{z\in \mc H: {\rm Re} \left<\tilde {\mc A} z, z\right>_{\mc H}= 0\right\}
\end{eqnarray}
has only the trivial solution, i.e. $z=0$; then by La Salle's invariance principle,
the system is strongly stable. In fact, ${\rm Re} \left<\tilde {\mc A} z, z\right>_{\mc H}=\left| \left( h_2 h_3 \varsigma \tilde BB_2 P_\varsigma  +  B_3 I\right)(z_2)_x(L) \right|^2=0.$  For letting $u=P_\varsigma (z_2)_x,$ $(z_2)_x=\tilde C \varsigma u - u_{xx},$
  $\tilde C \varsigma u(L)= u_{xx} (L)$ by the definition of $P_\varsigma$ in (\ref{def}), $u(L)=u_{xx}(L)=0.$ Thus,
$(z_2)_x(L)=[P_\varsigma (z_2)_x](L)\equiv 0$ by (\ref{set}).

Proving the strong stability of (\ref{Semigroupp-f}) reduces to showing that the following eigenvalue problem $\mc A z = \lambda z:$
\begin{eqnarray}
 \label{eig1}\left\{\begin{array}{ll}
 {\tilde A} w_{xxxx}+\gamma\beta \varsigma h_2 h_3 \tilde B^2 (J_\varsigma w_{x})_x +\lambda^2 w = 0, &\\
 w(0)=w_x(0)= w_{x}(L)= w_{xx}(L)=0, &\\
  {\tilde A} w_{xxx}(L)+\gamma\beta \varsigma h_2 h_3 \tilde B^2 J_\varsigma w_{x}(L)=(P_\varsigma w_x)(L)=0.
 \end{array}\right.
\end{eqnarray}
has only the trivial solution.
By using the definition of (\ref{MM}), i.e. $(J_\varsigma w_{x}) =(\varsigma \tilde C P_\varsigma w_x)-  w_x,$ we obtain that $(J_\varsigma w_{x})(L)=0$  since both terms $(P_\varsigma w_x)(L) $ and $w_x(L)$ are zero by (\ref{set}).

Let $\lambda=i\omega$ where $\omega\in\mathbb{R}.$ Then  (\ref{eig1}) reduces to
\begin{eqnarray}
 \left\{\begin{array}{ll}
 {\tilde A} w_{xxxx}+\gamma\beta \varsigma h_2 h_3 \tilde B^2 (J_\varsigma w_{x})_x -\omega^2 w = 0, &\\
w(0)=w_x(0)= w_{x}(L)= w_{xx}(L)=0, &\\
w_{xxx}(L)=J_\varsigma w_{x}(L)=(P_\varsigma w_x)(L)=0.
 \end{array}\right.
\end{eqnarray}
Note that the following integrals hold true.
\begin{eqnarray}
\begin{array}{ll}
 \int_0^L x w_{xxxx} \bar w_{xxx} dx=\frac{-1}{2}\int_0^L |w_{xxx}|^2dx, &\\
 \int_0^Lx  w \bar w_{xxx} dx= \int_0^L \frac{3}{2}\int_0^L |w_x|^2dx,&\\
 \int_0^L x(J_\varsigma w_x)_x\bar w_{xxx} dx= \int_0^L x((\tilde C \varsigma P_\varsigma -I)w_x)_x  \bar w_{xxx}dx&\\
\quad =  \int_0^L \tilde C\varsigma (P_\varsigma w_x)_x x\bar w_{xxx}dx +\frac{1}{2} \int_0^L |w_{xx}|^2dx.
\end{array}
\end{eqnarray}
Let $z=P_\varsigma w_x.$ Then $\tilde C\varsigma z -z_{xx}=w_x,$ and therefore
\begin{eqnarray}
\begin{array}{ll}
 \int_0^L \tilde C \varsigma (P_\varsigma w_x)_x x\bar w_{xxx} dx =\int_0^L \xi z_x x (\tilde C\varsigma \bar z_{xx}-\bar z_{xxxx})&\\
~~~=\frac{-1}{2}\int_0^L \left((\tilde C \varsigma)^2|z_x|^2 +\tilde C\varsigma |z_{xx}|^2\right) dx.
\end{array}
\end{eqnarray}
Multiplying  the equation (\ref{eig1}) by $x\bar w_{xxx}$ and integrate by parts and using the boundary conditions yields
\begin{eqnarray}
\nonumber \int_0^L \left[\tilde A |w_{xxx}|^2 +3 m |w_x|^2+ (\tilde C\varsigma)^2|z_x|^2 +\tilde C\varsigma |z_{xx}|^2\right] dx=0.
\end{eqnarray}
By using the overdetermined boundary conditions (\ref{eig1}) we obtain $w\equiv 0.  ~\square$

We state the following stability theorem:
\begin{thm}
   Then the solutions $\Phi$  for $t\in \mathbb{R}^+$ of the closed-loop system (\ref{Semigroupp-f})
     is exponentially stable in $\mc H.$
\end{thm}

\noindent {\bf Proof:} We prove the result by contradiction. Suppose that there exists  a sequence of real numbers $\beta_n\to \infty$ and a sequence of vectors $z_n=(w_n,v_n)\in {\rm Dom}(\mc A)$ with $\|z_n\|_{\mc H}=1$ such that
$ \|(i\xi_n  I-\mc A)z_n\|_{\mc H}\to 0, \quad  {\rm as}\quad n\to\infty, \quad {\rm i.e.}$
\begin{eqnarray}\left\{
\label{dal}\begin{array}{ll}
 i\xi_n w_n -v_n=f_n\to 0\quad {\rm in}\quad H^2_L(0,L)\\
 i\xi_n v_n +\frac{{\tilde A}}{m} (w_n)_{xxxx}\\
 \quad +\frac{\gamma\beta \varsigma h_2 h_3}{m} \tilde B^2 (J_\varsigma (w_n)_{x})_x=g_n\to 0\quad {\rm in}\quad \Ltwo(0,L).
 \end{array} \right.
\end{eqnarray}
By using the dissipation relationship (\ref{eq7}), we have
\begin{eqnarray}
\begin{array}{cc}
    i\xi_n \|w_n\|_{H^2_L(0,L)}^2 -\left<w_n,v_n\right>_{H^2_L(0,L)}&\\
    \quad\quad=\left<f_n,w_n\right>_{H^2_L(0,L)}=o(1),&\\
    i\xi_n \|v_n\|^2_{\Ltwo(0,L)} + <w_n, v_n>_{H^2_L(0,L)}&\\
   \quad\quad=\left<g_n,v_n\right>_{\Ltwo(0,L)}-d_n=o(1), \quad {\rm where }
\end{array}
\end{eqnarray}
\begin{eqnarray}
\begin{array}{ll}
\label{on}  d_n={\rm Re} \left<\left(
                                                 \begin{array}{c}
                                                   f_n \\
                                                   g_n \\
                                                 \end{array}
                                               \right),\left(
                                                 \begin{array}{c}
                                                   w_n \\
                                                   v_n \\
                                                 \end{array}
                                               \right)
\right>_{\mc H}&\\
 =\frac{k_1\gamma}{B_4}\left| \left( h_2 h_3 \varsigma \tilde BB_2 P_\varsigma  +  B_3 I\right)(v_n)_x(L)\right|^2.
 \end{array}
\end{eqnarray}
This implies that
\begin{eqnarray}
\label{dal1} \|w_n\|_{H^2_L(0,L)}^2-\|v_n\|^2_{\Ltwo(0,L)} =o(1).
\end{eqnarray}
Since $\|v_n\|_{\mc H}=1, $ and (\ref{dal}),(\ref{dal1}), we obtain
$ \|w_n\|_{H^2_L(0,L)}^2=\|v_n\|^2_{\Ltwo(0,L)} =\|\xi_n w_n\|_{\Ltwo(0,L)}^2=1/2.$
We need the following lemma to get a contradiction. The proof is provided in \cite{Ozer18c} due to the space limitation:

\begin{lem}\label{lemmai} Let $w_n\in {\rm Dom}(\mc A).$ Then, we have the following
\begin{eqnarray}
\begin{array}{ll}
\lim_{n\to \infty} \tilde A (w_n)_{xxx}(L) +\gamma\beta \varsigma h_2 h_3 \tilde B^2 D_x J_\varsigma (w_n)_x(L)=0,&\\
 \lim_{n\to \infty} \xi_n w_n(L)=\lim_{n\to \infty} (w_n)_{xx}(L)=0. \square
\end{array}&&
\end{eqnarray}
\end{lem}

Next, we simplify  (\ref{dal}) to get
 \begin{eqnarray}
\nonumber &&-\beta^2 _n w_n +\frac{{\tilde A}}{m} (w_n)_{xxxx}+\frac{\gamma\beta \varsigma h_2 h_3}{m} \tilde B^2 (J_\varsigma (w_n)_{x})_x\\
\label{dal3} &&~~ =i\beta_nf_n + g_n.
\end{eqnarray}

Let $q(x):=e^{x}.$  By taking the inner product of (\ref{dal3}) by $q (w_n)_x$ in $\Ltwo(0,L)$ to get
 \begin{eqnarray}
\label{dal4} \begin{array}{ll}
 \left<i\beta_nf_n + g_n, q (w_n)_x\right>=\left<-\beta^2 _n w_n +\frac{\tilde A}{m}(w_n)_{xxxx} \right.&\\
 \left.\quad\quad+\frac{\gamma\beta \varsigma h_2 h_3 \tilde B^2 (J_\varsigma (w_n)_{x})_x}{m}, q  (w_n)_x\right>\to 0
\end{array}
\end{eqnarray}
since there exists constants $D_1, D_2>0,$
\begin{eqnarray}
\begin{array}{ll}
 \left<g_n,q (w_n)_x\right>_{\Ltwo(0,L)}\le D_1 \|g_n\|_{\Ltwo(0,L)} \|w_n\|_{H^2_L(0,L)}\to 0,&\\
|\left<i\beta_nf_n ,q (w_n)_x\right>_{\Ltwo(0,L)}| \le & \\
 D_2\left(\|f_n\|_{H^2_L(0,L)}\|\beta_n w_n\|_{\Ltwo(0,L)} + |f_n(1)\beta_n w_n(1)|\right)\to 0
\end{array}
\end{eqnarray}
where we used Lemma \ref{lemmai}. By integration by parts,
\begin{eqnarray}
\label{d1}\begin{array}{ll}
 {\rm Re}\left<-m\beta_n^2  w_n, q (w_n)_x\right>_{H}=-\frac{m e}{2} |\beta_n w_n(1)|^2 &\\
  \quad +\frac{m}{2}\int_0^1 e^x|\beta_n w_n|^2dx,
\end{array}
\end{eqnarray}
\begin{eqnarray}
\label{d2}\begin{array}{ll}
  {\rm Re}\left< {\tilde A} (w_n)_{xxxx}+\gamma\beta \varsigma h_2 h_3 \tilde B^2 (J_\varsigma (w_n)_{x})_x,q (w_n)_x\right>&\\
 = {\rm Re}[{\tilde A} (w_n)_{xxx}+\gamma\beta \varsigma h_2 h_3 \tilde B^2 (J_\varsigma (w_n)_{x}&\\
 \left.\quad  -{\tilde A} e^x (w_n)_{xx} (w_n)_x\right](L)  - \frac{e\tilde A}{2} |w_n(L)|^2&\\
  \quad + \int_0^L \left\{\frac{3\tilde A}{2}e^x |(w_n)_{xx}|^2  + {\rm Re}[({\tilde A} e^x (w_n)_{xx}\right.&\\
\left. ~~~ +\gamma\beta \varsigma h_2 h_3 \tilde B^2 e^x(J_\varsigma (w_n)_{x})(\bar w_n)_x]\right\}dx.
\end{array}
\end{eqnarray}
The boundary terms converge to zero due to Lemma \ref{lemmai}, and since $\|w_n\|_V<\infty. $ Therefore $\|w_n\|_{\Ltwo(0,L)}=o(1),$ and
\begin{eqnarray}
\begin{array}{ll}
\|w_n\|_{H^1_L(0,L)}\le \sqrt{\|w_n\|_{\Ltwo(0,L)}\|w_n\|_{H^2_L(0,L)}}=o(1),&\\
{\rm and}\quad \int_0^L \tilde A e^x (w_n)_{xx}({\bar w}_n)_x dx=o(1)
\end{array}
\end{eqnarray}
Using (\ref{d1}) and (\ref{d2}) in (\ref{dal4})  we get $\|w_n\|_{H^2_L(0,L)}=o(1)$ contradicting  with $\|z_n\|_{\mc H}=1.$

\section{Stable Approximations \& Simulations}
\label{simul}
The aim of this section is to present a sample numerical experiment in order to show that  the stabilizing boundary controller (\ref{feed})  can be designed numerically. Since our model (\ref{abstractMM}) is strongly coupled,  it requires a more careful treatment for the high frequency modes which may cause spill-overs. The widely-used approximations, i.e. the standard Galerkin-based Finite Element or Finite Difference,  fail to provide reliable results  for boundary control problems \cite{B-I-W}. The filtering technique for Finite Differences has been recently developed to avoid artificial high-frequency solutions causing  instabilities in the approximated solutions. This is achieved by adding extra distributed damping terms to the equations or boundary conditions, as in \cite{Leon,Roventa,T-Z}.

We consider a three-layer smart beam with  length  $L=1 {\rm  m},$ and thicknesses of each layer $h_1,h_3=0.1{\rm m},$  $h_2= 0.01{\rm m}.$  The  material constants are chosen $\rho_1,\rho_3=7600$ kg/{m$^3$}, $\rho_2= 5000$ kg/{m$^3$},  $\alpha_1,\alpha_3=1.4\times 10^7$ N/{m{$^2$}},  $\alpha_2=10^5$ N/{m{$^2$}}, $\gamma=10^{-3}$ C/{m$^2$}, $\beta= 10^{6} {\rm m/F},$   $G_2=100$ GN/${\rm m}^2$. We consider the simulation for $T<5,$ and initial data $w(x,0)= \dot w(x,0)=10^{-4}\sum\limits_{i=2}^4 e^{\left(\frac{x-i*L/5}{0.2L}\right)^2}$.  We also non-dimensionalize the time variable  $t=A_1 t^*$  with $t^* \mapsto t$ where $A_1=L\sqrt{\frac{\rho}{\tilde A}}\sim 0.82.$   Now consider the  discretization of the interval $[0,L]$ with the fictitious points $x_{-1}$ and $x_{N+1}$ with $N=60:$
\begin{eqnarray}
\begin{array}{ll}x_{-1}<0=x_0< x_1 <x_2, \ldots < x_N = L < x_{N+1}, &\\
 x_i=i\cdot dx, \quad i=-1,0,1,\ldots, N, N+1, ~~dx=\frac{L}{N+1}.
 \end{array}
\end{eqnarray}
Henceforth, to simplify the notation, we use $z(x_{i})=z_i.$ We adopt the semi-discrete scheme in Finite Differences to simulate the effects of the stabilizing controller. The following are the second order finite difference approximations for different order derivatives:
\begin{eqnarray}
\nonumber \begin{array}{ll}
z_{x}=\frac{z_{i+1}-z_{i-1}}{2dx}, \quad {\rm or},~~ z_x=\frac{3z_i-4z_{i-1}+z_{i-2}}{2dx},&\\
 z_{xx}=\frac{z_{i+1}-2z_i+z_{i-1}}{dx^2},&\\
 z_{xxx}=\frac{z_{i+2}-2z_{i+1}+2z_{i-1}-z_{i-2}}{2dx^3}&\\
 z_{xxxx}=\frac{z_{i+2}-4z_{i+1}+6z_{i}-4z_{i-1}+z_{i-2}}{dx^4}.
\end{array}
\end{eqnarray}
  The numerical viscosity terms $-\dot w_{xx}$ and $-\dot \phi^2_{xx}$ are added to the $w$ and $\phi^2-$equations in (\ref{abstractMM-damped}) , respectively. The discretization of (\ref{abstractMM-damped}) is
\begin{eqnarray}
\nonumber
\begin{array}{ll}
   \ddot w  +\frac{w_{i+2}-4w_{i+1}+6w_i-4w_{i-1}+w_{i-2}}{dx^4}&\\
  ~ -\kappa\frac{\dot w_{i+1}-2\dot w_i+\dot w_{i-1}}{2dx^2}+\frac{\beta  \gamma  \varsigma  h_2 h_3 L^3 \tilde B}{\tilde A}\frac{\phi^2_{i+1}-\phi^2_{i-1}}{2dx}=0,&\\
   -\frac{\dot \phi^2_{i+1}-2\dot \phi^2_{i}+\dot \phi^2_{i-1}}{L^2dx^2}-\frac{\phi^2_{i+1}-2\phi^2_{i}+\phi^2_{i-1}}{L^2dx^2}&\\
   ~~+\varsigma \tilde C \phi^2_i +\frac{\tilde B}{L^3}\frac{w_{i+2}-2w_{i+1}+2w_{i-1}-w_{i-2}}{2dx^3}=0, &\\
   \qquad i=1,\ldots, N-1,&\\
 \phi^2_0=w_0=0, w_{-1}=w_1, &\\
 \frac{3\phi^2_N-4\phi^2_{N-1}+\phi^2_{N-2}}{2dx} -\left(\frac{\gamma \tilde B B_3}{L^2 B_4 \tilde A}-\frac{B_2}{\beta B_4}\right)k_1V(t)=0,&\\
    \frac{w_{N+1}-2w_N+w_{N-1}}{dx^2} -\frac{\gamma B_3}{B_4 \tilde A}k_1V(t)=0,&\\
 \tilde A\frac{2w_{N+1}-5w_N+2w_{N-1} + 4w_{N-2}-4w_{N-3}}{2dx^3}&\\
 \quad +\beta  \gamma  \varsigma  h_2 h_3 L^3 \tilde B \phi^2_N=0&\\
\end{array}
\end{eqnarray}
where $\kappa=\frac{dx}{5}<dx,$ and the Voltage controller $V(t)$ is designed as the following (with the choice of $k_1=10^8$)
\begin{eqnarray}
\nonumber
  \begin{array}{ll}
    V(t) =\frac{ \varsigma  \tilde C+\tilde B}{\varsigma  \tilde C} \frac{3\dot w_N-4\dot w_{N-1}+\dot w_{N-2}}{2dx} -\frac{1}{\varsigma \tilde B\tilde C }\dot\phi^2_N  \end{array}.
\end{eqnarray}


\begin{figure}[htb!]
    \centering
    \vspace{-0.02in}
    {{\includegraphics[width=4.3cm]{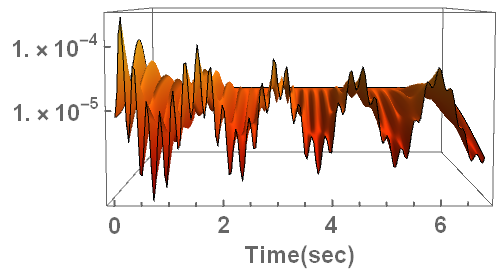} }}%
            {{\includegraphics[width=4.3cm]{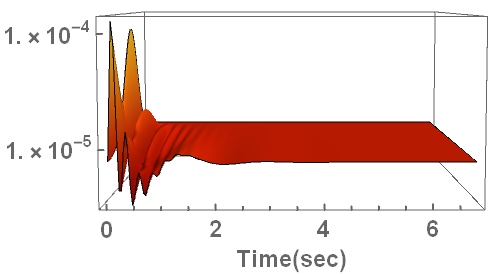} }}%
        \vspace{-0.1in}
    \caption{Rapid decay of the bending  $w(x,t)$ in a few seconds (real time) after the controller applies. }%
    \label{w}%
\end{figure}
\begin{figure}[htb!]
    \centering
    \vspace{-0.02in}
    {{\includegraphics[width=4.4cm]{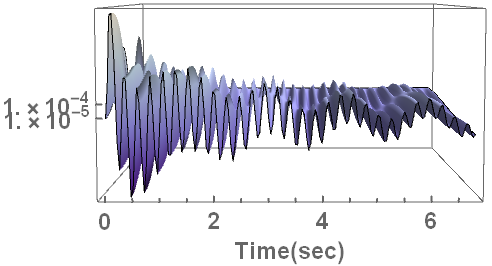} }}%
            {{\includegraphics[width=4.3cm]{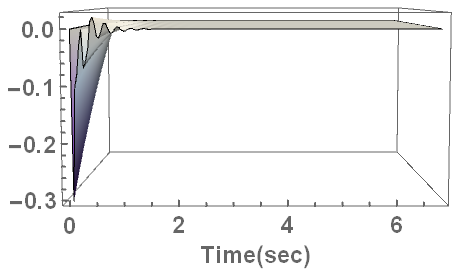} }}%
        \vspace{-0.1in}
    \caption{Rapid decay of the shear  $\phi^2(x,t)$ in a few seconds (real time) after the controller applies. }%
    \label{phi}%
\end{figure}

\begin{figure}[htb!]
    \centering
    \vspace{-0.02in}
    \includegraphics[width=4.2cm]{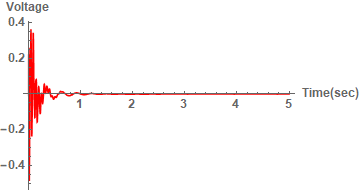} 
           \includegraphics[width=4.2cm]{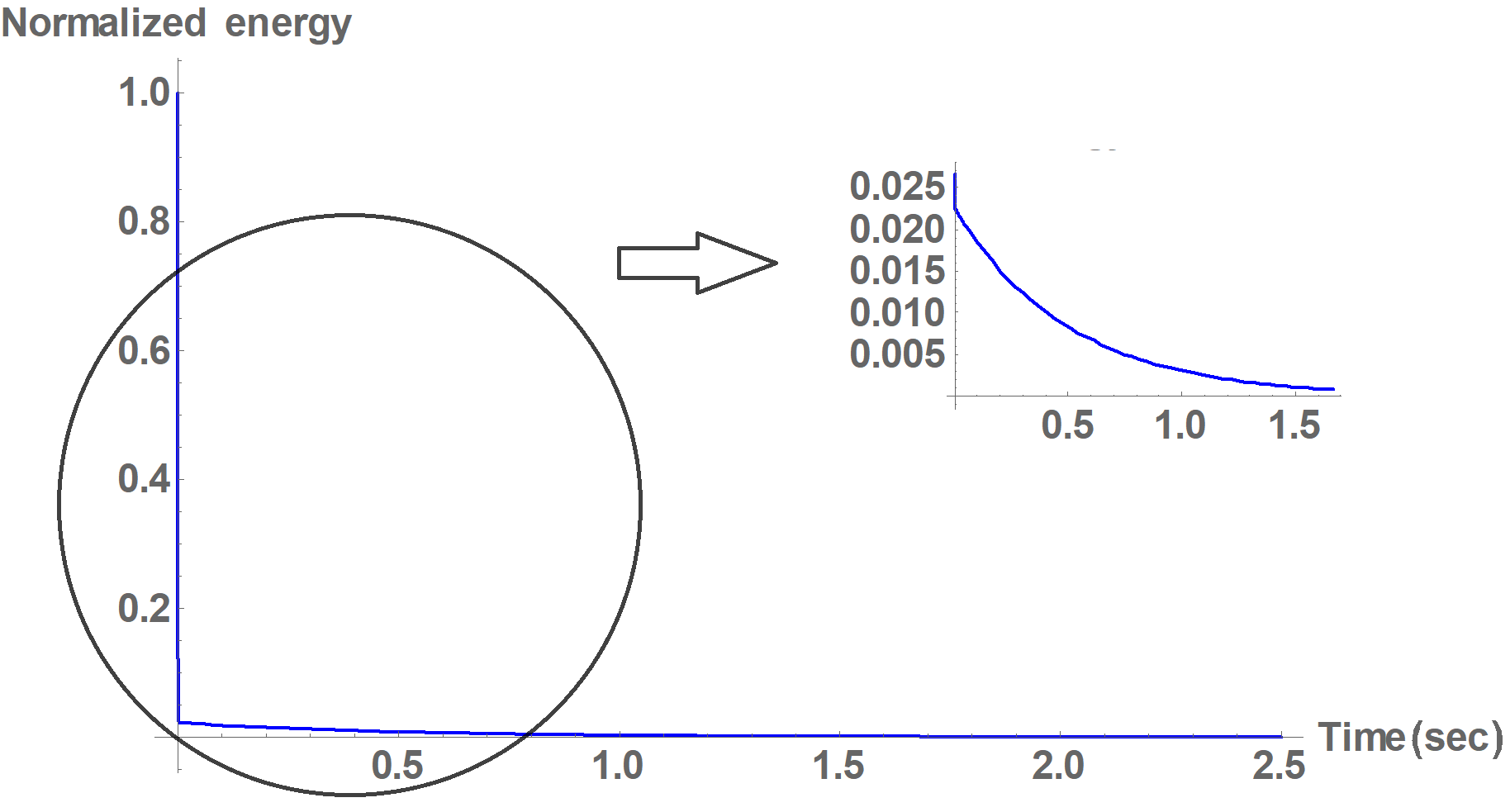}
        \vspace{-0.1in}
    \caption{Voltage $V(t)$  and normalized energy $E(t)$ distributions for the first  few seconds. }%
    \label{voltage}%
\end{figure}

 The simulations in Figures \ref{w} and \ref{phi} show that the $\phi^2$ and $w$ solutions both decay to zero fast enough. In fact, $\phi^2$ solution destabilizes in the beginning (the picking phenomenon in Fig. \ref{voltage}) but then it decays to zero faster than the bending solution. These results can be tuned up by using an improved scheme after a careful stability analysis is performed.

Note the necessity of the controller $P_\varsigma \dot w_x(L)$ in (\ref{feed}) to prove the strong stability result in Theorem \ref{strong}. It is an open problem to analytically prove the same result without $P_\varsigma \dot w_x(L).$
 In fact, further numerical investigation is the subject of \cite{Ozer18c} where the impact of the non-classical feedback controller $V_1(t)=-k_2(P_\varsigma \dot w_x)(L)$   over the classical one $V_2(t)=-k_1\dot w_x(L)$   is shown to be crucial (different feedback gains for each).

Models incorporating the nonlinear elasticity theory are also derived by a consistent variational approach, and  the filtering technique is applied in \cite{Ozer18}. The reader should refer to promising numerical results in \cite{Ozer18} with the choice of various nonlinear stabilizing feedback controllers.  The results of this paper   will be a basis for  functional and numerical analyses for the nonlinear beam models in  \cite{K-O}. Developing stable Finite Difference schemes and the adoption of the mixed-Finite Element method for both linear and nonlinear models  are the progressing works \cite{Ozer18c}.


\end{document}